\documentclass[pra,preprint]{revtex4}
\usepackage{epsfig}
\usepackage{amsmath}
\usepackage{epsfig}

\begin{document}
\title
{Families of curves orthogonal to the lines $y=mx-2m-m^3$}
	\author{Zafar Ahmed$^{1,2}$ and Pallavi S. Telkar$^3$}
	\affiliation{$^1$Nuclear Physics Division, Bhabha Atomic Research Centre, Mumbai 400 085, 
	 $^4$Department of Mathematics, Institute of Chemical Technology, Mumbai, India}
	\email{1,2:zahmed@barc.gov.in, 3: telkarpallave@gmail.com}
\date{\today}
\begin{abstract}
The family of lines $y=mx-2m-m^3$, are well known to be normal to the parabola $y^2=4x$. However, this family of lines is normal to a family of curves of which this parabola is just one member.
Here, by solving an interesting first order and third degree ODE, we bring out these curves. The resulting one set of curves are ``parabola-like" but non-standard ones and the other family is not even ``parabola like".
\end{abstract}
\maketitle
Orthogonal pairs of family of curves represent a co-ordinate system.  They may also represent  an electrostatic system  comprising of lines of force and equi-potentials. 
Textbooks $[1,2]$ discuss the method and exercises of finding orthogonal family of curves to a given one-parameter  family of curves. First,  one sets up a first order  ordinary differential equation (ODE) by eliminating the free parameter of the given family of curves, by differentiation. Next, in this ODE(c) one changes $y'$ to $\frac{-1}{y'}$ to set up a new ODE(o) to solve and get the required family of their orthogonal curves. Here c denotes the curve and o denotes the orthogonal curve/trajectory.

The orthogonal trajectories to the curves $x^2-y^2=C_1$ is found by differentiating this equation and setting up the ODE(c) as $y'=\frac{x}{y}$, the ODE(o)  is obtained by changing $y'$ to $-1/y'$, we get ODE(o) as $y'=-\frac{y}{x}$, solving this, we get $xy=C_2$, where $C_2$ is  the constant/parameter of integration.  The example just discussed represents a hyperbolic co-ordinate system and  an electrostatic quadrupole as well. More involved orthogonal pairs of curves are discussed in the textbooks [1,2]. 

Let us find orthogonal trajectories to the family of lines $y=mx$, where $m$ is a real parameter. Differentiating $\frac{y}{x}=m$, we eliminate $m$ and get the ODE(c) as $y'=\frac{y}{x}$. Next we change $y'$ to $-1/y'$ to get ODE(o) as $y'=-\frac{x}{y}$, solving this we find orthogonal trajectories as $x^2+y^2=r^2$, where $r$ is a real parameter. This pair of curves represent the polar co-ordinate system (when, $m=\tan \theta$) and the electric mono-pole (charge) where the lines are the lines of force and the circles are equi-potentials.

Similarly, we can find the orthogonal trajectories to the family of lines $y=m(x+1)$, The ODE(c) for this family is obtained by eliminating $m$ by differentiation $\frac{y}{x+1}=m$, we get $(x+1)y'-y=0$. By changing $y'$ to $-1/y'$, we get ODE(o) as $(x+1)+yy'=0$, solving this  we get $(x+1)^2+y^2=r^2$, where $r$ is a real parameter. For  any real value of $m$ the lines are normal (orthogonal) to the circles.

We would like to remark that finding orthogonal trajectories to  other one-parameter families of lines: $y=mx+f(m)$, where $f(m)$ becomes a more involved function of $m$, remain elusive in the literature. The well known  one-parameter family of lines [3]
\begin{equation}
y=mx-2m-m^3.
\end{equation}
are known to be normal to the simple parabola $y^2=4x$ [3]. However, other curves  which are orthogonal to  the lines (1)  are not known. In this article, we propose to find them. Let us differentiate (1) w.r.t. $x$, we get $y'=m$, so the required ODE(c) is
\begin{equation}
y=y'x-2y'-{(y')}^3.
\end{equation}
The ODE(o) for the family of  curves orthogonal to (1) can be obtained by replacing $y'$ by $-1/y'$.
\begin{equation}
y(y')^3 = (y')^2 (2-x)+1 
\end{equation}
Denoting $y'$ as $p$, we get,
\begin{equation}
yp^3 = p^2 (2-x)+1
\end{equation}
differentiating w.r.t $y$, we get,
\begin{equation}
3p^2y \frac{dp}{dy}+p^3 = -2px \frac{dp}{dy}-p^2 \frac{dx}{dy} +4p \frac{dp}{dy}
\end{equation}
Substituting $x$ from Eq.$(4)$ and $\frac{dx}{dy}=\frac{1}{p}$  in Eq.$(5)$, we get an ODE in only two variables $y, p$:
\begin{equation}
(p^3+p)dy + (yp^2+\frac{2}{p})dp =0
\end{equation}
Comparing (6) with $Mdy+Ndp=0$, $\frac{\partial M}{\partial p}=3p^2+1, \frac{\partial N}{dy}=p^2$.
Since the differential Eq.$(6)$ is not exact [1,2], it can be made so by multiplying with the integrating factor $\mu(p)$, where $\mu(p) $ [1,2] is
\begin{equation}
\mu(p) =\exp \int \left(\frac{\frac{\partial N}{dy}-\frac{\partial M}{dp}}{M} \right) dp= \exp{\left(\int\frac{-2p^2-1}{p^3 +p}dp\right)} = \frac{1}{p \sqrt{1+p^2}}
\end{equation}
Multiplying Eq.(6) by $\mu(p)$ in Eq.(7), we get the exact ODE,
\begin{equation}
\sqrt{1+p^2}dy + \left(py+\frac{2}{p^2}\right)\frac{1}{\sqrt{1+p^2}}dp = 0.
\end{equation}
The solution of Eq.(8) can be written as 
\begin{equation}
\int \sqrt{1+p^2}~ dy ~ \text{[treat $p$ as constant]} + \int \frac{2}{p^2} \frac{1}{\sqrt{1+p^2}}~ dp~ \text{[terms of $N$ not containing $y$]}=C
\end{equation}
Integrating (9), we get
\begin{equation}
y\sqrt{1+p^2}-\frac{2 \sqrt{1+p^2}}{p} = C
\end{equation}
Thus,
\begin{equation}
y = \frac{2}{p}+\frac{C}{\sqrt{1+p^2}} 
\end{equation} 
where C is a constant of integration. Inserting (11) in (4), we get
\begin{equation}
x = \frac{1}{p^2} - \frac{C p}{\sqrt{1+p^2}}
\end{equation}
Eq.(11) and Eq.(12) are the parametric solution of Eq.(3), where eventually $p$ acts merely as a real  parameter. By choosing $t = 1/p$, we get more convenient parametrization of family of curves which are orthogonal to the family of lines in Eq.(1). We write,

\begin{equation}
x = t^2 - \frac{C}{\sqrt{1+t^2}} ,\quad  y = 2t+\frac{Ct}{\sqrt{1+t^2}}~.
\end{equation}

\begin{figure}[ht]
	\centering
	\includegraphics[width=8cm,height=9cm,scale=1]{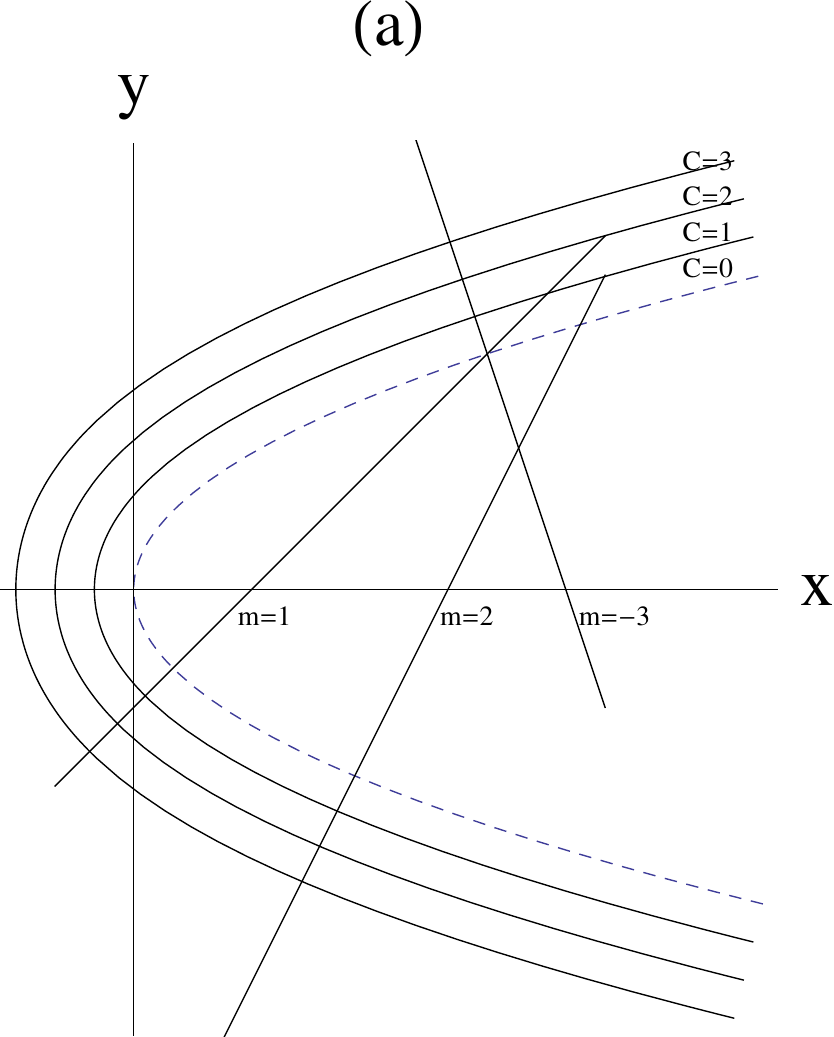}
	\hskip.2cm
		\includegraphics[width=8cm,height=9cm,scale=1]{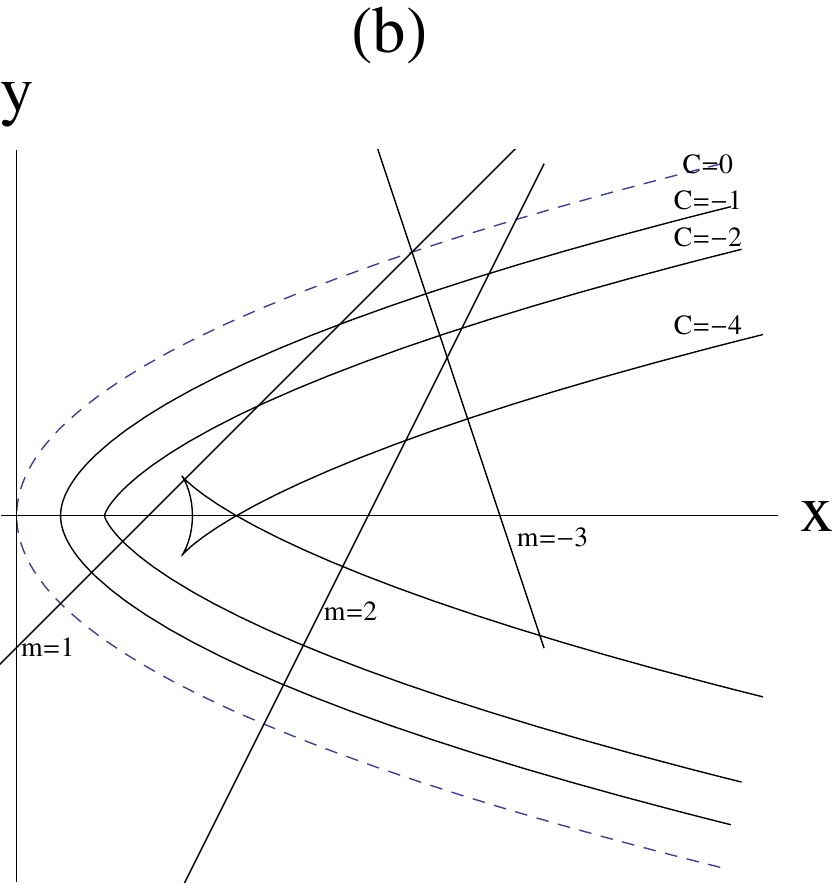}	
	\caption{Depicting 8 members of the family of curves (13) orthogonal to the family of lines (1)
	$(m=1,2,-3)$. Only the dashed curve is parabola others are all non-parabolic.} 
		
\end{figure}

Where the most popular curve, parabola $(y^2 = 4x)$ emerges when $C= 0$. Eventually, any member of the family of lines (1) is orthogonal to any of the curve (13) only at one point. The line may cut the curve non-orthogonally at some other point.

Now the  question is what the non parabolic family of curves would look like.
In Fig. 1, we give a graphical depiction of the parabolic ($C=0$) and the non-parabolic $(C \ne 0$) curves which are always orthogonal to the lines (1) for any real value of m.

Finally, we conclude that the family of orthogonal curve(s) (13) to the family of lines (1) is the  parabola $y^2=4x$  when $(C=0)$ else for $C\ne 0$ they are new and non-parabola even if they look like a parabola. However, surprisingly, for large negative values of $C$  they are not even parabola-like, see the  inner most curve in Fig. 1(b) for $C=-4$. Eventually, for $C \ne 0$ the orthogonal curves(trajectories) (13) are new, non-standard and non-parabolic,  It may be remarked that  a parabola is always represented as  a (rational function of $x,y$)  quardratic: $(y-m_1x-c_1)^2=A(y-m_2x-c_2)$, where $m_1 \ne m_2.$ If we eliminate $t$ in (13), for $C\ne 0$ we do not get a rational relation between
$x$ and $y$.

\vskip 1 cm
\noindent
[1] See for example, G. F. Simmons and S. G. Krantz, `Differential Equations' (McGraw-Hill, New York) 2007.

\noindent
[2] See for example, M. L. Krasnov, A.I. Kiselyov, G.I. Makarenko; `A Book of Problems in Ordinary Differential Equations' (MIR, Moscow) 1978.

\noindent
[3] See for example, S. L. Loney, `Elements of Co-ordinate geometry,' Part 1.
\end{document}